\documentclass{elsart}
\usepackage{euscript,epsfig}
\usepackage{amssymb}
\usepackage{amsmath}
\usepackage{graphicx}
        \newcommand\beqn{\begin{eqnarray*}}
        \newcommand\eeqn{\end{eqnarray*}}
        \newcommand\beqa{\begin{eqnarray}}
        \newcommand\eeqa{\end{eqnarray}}
        \newcommand\br{\begin{array}}
        \newcommand\er{\end{array}}
        \newcommand\bes{\begin{equation*}}
        \newcommand\ees{\end{equation*}}
        \newcommand\be{\begin{eqnarray}}
        \newcommand\ee{\end{eqnarray}}
        \newcommand\efig{\end{figure}}
        \newcommand\btab{\begin{table}}
        \newcommand\etab{\end{table}}
        \newcommand\btabl{\begin{tabular}}
        \newcommand\etabl{\end{tabular}}
        \newcommand\bm{\begin{math}}
        \newcommand\ema{\end{math}}
        \newcommand\bi{\begin{itemize}}
        \newcommand\ei{\end{itemize}}
        \newcommand\bd{\begin{description}}
        \newcommand\ed{\end{description}}
        
        \newcommand{\bpr}{\begin{proof}}
        \newcommand{\epr}{\end{proof}}
        \newtheorem{lemma}{\bf Lemma}
        \newtheorem{theorem}{\bf Theorem}

         \newcommand{ \vect }[2]{
                    \left( \begin{array}{c}
                           #1 \\
                           #2
                            \end{array} \right)}

        \newcommand\Bc{\Eu {B}_c (\ell^2, \ell^2)}

         \newcommand\hV{\hat{V}}
           \newcommand\hU{\hat{U}}
             \newcommand\hM{\hat{M}}
               \newcommand\hN{\hat{N}}

        \newcommand\Eu{\EuScript}
        \newcommand\mrm{\mathrm}

\newcommand\EuBt{\Eu B(\ell^2, \; \ell^2 \times \ell^2)}

\date{}
\begin{document}

\begin{frontmatter}
\title{On Robustness in the Gap Metric and Coprime Factor Uncertainty for LTV Systems}
\author{Seddik M. Djouadi}
\thanks{This work was partially supported by a Ralph E. Powe  Junior Enhancement award.}
\address{Department of Electrical Engineering and Computer Science,
\\ University of Tennessee, Knoxville, TN 37996-2250,
\\ djouadi@eecs.utk.edu}
\begin{abstract}
In this paper, we study the problem of robust stabilization for linear time-varying (LTV) systems subject
to time-varying normalized coprime factor uncertainty. Operator theoretic results which generalize similar
results known to hold for linear time-invariant (infinite-dimensional) systems are developed. In particular,
we compute an upper bound for the maximal achievable stability margin under TV normalized coprime
factor uncertainty in terms of the norm of an operator with a time-varying Hankel structure.
We point to a necessary and sufficient condition which guarantees compactness of the TV Hankel operator,
and in which case singular values and vectors can be used to compute the time-varying stability margin and
TV controller. A connection between robust stabilization for LTV systems and an Operator Corona Theorem
is also pointed out.
\end{abstract}
\end{frontmatter}
\section*{Definitions and Notation}
\label{notation}
\begin{itemize}
\item  $\Eu {B} (E, F)$ denotes the
space of bounded linear operators from a Banach space $E$ to a Banach
space $F$, endowed with the operator norm
\beqn
\| A\| :=\sup_{x\in E, \; \|x\| \leq 1} \|Ax\|, \;\; A\; \in\;\Eu {B} (E, F)
\eeqn
\item $\ell^2$ denotes the usual Hilbert space of square summable sequences
with the standard
norm
\beqn
\| x \|_2^2 := \sum_{j=0}^\infty |x_j|^2,
\;\; x := \bigl(x_0, x_1, x_2, \cdots \bigl) \in \ell^2
\eeqn
\item $P_{k}$ the usual truncation operator for some integer $k$, which sets all
outputs after time $k$ to zero.
\item An operator $A \in \Eu {B} (E, F)$ is said to be causal if
it satisfies the operator equation:
\beqn
P_{k }A P_{k } = P_{k } A, \; \forall k \;\; {\rm positive\;\;integers}
\eeqn
\item $tr(\cdot)$ denotes the trace of its argument.
\end{itemize}
The subscript ``$_c$'' denotes the restriction of a subspace of operators to its intersection
with causal (see \cite{saeks,feintuch} for the definition) operators. ``$\oplus$'' denotes for the
direct sum of two spaces. ``$^\star$'' stands for the adjoint of an operator or the dual space of
a Banach space depending on the context
\cite{douglas,luen}.
\section{Introduction}
The gap metric was introduced to study stability robustness of feedback systems. It induces the weakest topology in which
feedback stability and performance (as measured by in terms of a closed-loop induced norm) are robust \cite{sakkary,georgiou1,smith,GS2,zhu,vinni}.
In \cite{georgiou1} Georgiou showed the relationship
between the gap metric and a particular two-block $H^\infty$ problem. In \cite{smith}, the authors showed that feedback
optimization in the gap metric is equivalent to feedback optimization with respect to normalized factor perturbations. They
computed the largest possible uncertainty radius such that robust stability is preserved. Extensions to time-varying systems
have been proposed in  \cite{foiascdc,georgiou} where a geometric framework for robust stabilization of infinite-dimensional
time-varying systems was developed. The uncertainty was described in terms of its graph and measured in the gap metric.
Several results on the gap metric and the gap topology were established, in particular, the concept of a graphable subspace
was introduced. In \cite{feintuch}, some of the results obtained in \cite{georgiou1} were generalized, in particular, the gap
metric for time-varying systems was generalized to a two-block time varying optimization analogous to the two-block $H^\infty$-
optimization proposed in \cite{georgiou1}. This was achieved by introducing a metric which is the supremum of the sequence of gaps
between the plants measured at every instant of time. The latter reduces to the standard gap metric for linear time-invariant
(LTI) systems. \\ \\
In \cite{fein,feintuch,avraham} using the time-varying gap metric it is shown that the ball of uncertainty in
the time-varying gap metric of a given radius is equal to the ball of uncertainty of the same radius defined by perturbations
of a normalized right coprime fraction, provided the radius is smaller than a certain quantity. In \cite{avraham} lower
and upper bounds are derived for computing the maximal stability margin (quantified in the TV gap metric), in terms of
using coprime factorizations. These bounds are equal to the maximal stability margin for LTI systems. In \cite{djouadi},
the authors showed that the time varying (TV) directed gap reduces to the computation of an operator with a TV Hankel plus Toeplitz structure. Computation of the norm of such an operator can be carried out using an iterative scheme known to hold
for standard two-block $H^\infty$ problems \cite{djouadi}. The minimization in the TV directed gap formula was shown to be a minimum using
duality theory. \\ \\
In this paper, we use the equivalence between uncertainty quantified in TV gap metric balls and coprime
factor uncertainty for LTV systems. We study the problem of robust stabilization for time-varying normalized coprime factor
perturbations and obtain operator theoretic results, which generalize similar results in \cite{glover,smith,GS2} known to hold
for LTI systems. In particular, we compute an upper bound on the maximal achievable stability margin under TV normalized
coprime factor uncertainty in terms of the norm of a TV Hankel operator. The upper bound reduces to the maximal stability margin for LTI systems.
We point to a necessary and sufficient condition which guarantees compactness of the TV Hankel operator,
and in which case singular values and vectors can be used to compute the TV optimal stability margin and TV controller.
Therefore, generalizing similar results obtained in \cite{GS2} for LTI systems to TV ones. The technique developed to compute the upper
bound also applies to the lower bound on the maximal achievable stability margin in the TV case. To this end, it suffices
to restrict the relevant operators to a particular subspace. Computing an upper and lower bound allows the estimation of the optimal stability margin
for LTV systems. We point to an Operator Corona Theorem connection which gives a necessary
and sufficient condition for the existence of a robustly stabilizing controller.
\\
\\
The rest of the paper is organized as follows. In section \ref{s1} the gap metric is introduced. In section
\ref{s2} the relation between the TV gap metric and coprime factorization is discussed,
and computations in terms of operator theory are developed. We conclude with a summary of our contribution in section \ref{s4}.
\section{The Time-Varying Gap Metric} \label{s1}
LTV systems may be regarded as causal linear (possibly unbounded) operators acting on $\ell^2$ as multiplication operators.
To each plant $P$ we associate its domain
\beqa
D(P) = \{u\in \ell^2: Px \in \ell^2\} \label{1}
\eeqa
An LTV plant $P$ has a right coprime factorization if there exist operators $M$ and $N$ both in $\Eu {B}_c (\ell^2, \ell^2)$,
such that, $P=NM^{-1}$, and a left coprime factorization if there exist $\hat{M}, \hat{N} \in \Eu {B}_c (\ell^2, \ell^2)$,
such that $P=\hat{M}^{-1} \hat{N}$. In addition there exist $X$, $Y$, $\hat{X}$, $\hat{Y}  \in \Eu {B}_c (\ell^2, \ell^2)$,
such that,
\beqa
XN + YM = I, \;\;\;\; \hat{N}\hat{X} + \hat{M}\hat{Y}  = I
\label{4}
\eeqa
Such factorizations exist if and only if $P$ is stabilizable \cite{dale}.
There exist causal bounded linear operators $U, V, \hat{U},$ and $\hat{V}$ such that \cite{dale}
\beqa
\left( \br{cc} \hV & -\hU  \\ -\hN
& \hM \er \right)
\left( \br{cc} M & U  \\ N
& V \er \right)
= \left( \br{cc} M & U  \\ N
& V \er \right)  \left( \br{cc} \hV & -\hU  \\ -\hN
& \hM \er \right)
=  \left( \br{cc} \;I\; & \;0\;  \\ \;0\;
& \;I\; \er \right)
\label{6}
\eeqa
All $P$ stabilizing LTV controllers $C$ can be parameterized as \cite{dale}
\beqa
C = (U+MQ)(V+NQ)^{-1} = (\hV +Q\hN)^{-1} (\hU + Q\hM), \; Q \in \Eu {B}_c (\ell^2, \ell^2)
\label{7}
\eeqa
Following \cite{smith,feintuch} we are interested in normalized coprime factorizations. That is,
$\vect{M}{N}$ is an isometry from $\ell^2$ into $\ell^2 \oplus \ell^2$, in this case,
$M^\star M+N^\star N=I$ as is in the LTI case, where $M^\star$, $N^\star$ are the adjoint operators of $M$ and
$N$, respectively \cite{feintuch}.
\\
\\
Suppose that two LTV plants $G_1$ and $G_2$ have normalized right coprime factorizations $\vect{M_1}{N_1}$
and $\vect{M_2}{N_2}$, respectively.
Following \cite{fein,feintuch} Let $\Pi_{1n}$ denote the orthogonal projection on the range of
$\vect{M_1}{N_1}(I-P_n)$, $\Pi_{2n}$ is defined similarly. Define
\beqa
\overrightarrow{\delta_n}(G_1, G_2) = \left\| \left\{ \left( \br{cc} I-P_n & -0  \\ 0
& I-P_n \er \right) - \Pi_{2n}\right\} \Pi_{1n} \right\|
\label{10}
\eeqa
and
\beqa
\delta_n(G_1, G_2) = \| \Pi_{1n}- \Pi_{2n}\| = \max \bigl(
\overrightarrow{\delta_n}(G_1, G_2),
\overrightarrow{\delta_n}(G_2, G_1) \bigl)
\label{11}
\eeqa
The directed time varying gap between $G_1$ and $G_2$ is then defined as
\beqa
\overrightarrow{\alpha}(G_1, G_2) := \sup_{n\geq 0} \overrightarrow{\delta_n}(G_1, G_2)
\label{14}
\eeqa
and the time varying gap \cite{feintuch}
\beqa
\alpha(G_1, G_2) = \max\Bigl( \overrightarrow{\alpha}(G_1, G_2),
\overrightarrow{\alpha}(G_2, G_1)\Bigl)
\label{15}
\eeqa
The function $\alpha$ is a metric and for time-invariant systems reduces to the standard gap
metric $\delta$ \cite{fein}. 
\\\\
Let $B(P, r)$ denote the set of stabilizable LTV systems $P_1$ such that
\beqa
\alpha(P,\; P_1) \; < \; r
\eeqa
Denote by $B_s(P, r)$ the set of all $P$ with right coprime factorization $N M^{-1}$ for
which
\beqa
\left\| \vect{M}{M}- \vect{M_1}{N_1}\right\| < r
\eeqa
and such that $M$ is invertible in the algebra of LTV systems. Then the following result in \cite{feintuch} (Lemma 5.3. p. 219) holds.
\begin{theorem}(\cite{feintuch}, p. 219)   \label{t11}
\beqa
B(P, r) = B_s(P, r)
\eeqa
for $r < \frac{1}{\inf_{Q\in \Eu {B}_c (\ell^2, \ell^2)} \|[\hat{V}+Q\hat{N}, \; -(\hat{U}+Q\hat{M})]\||}$.
\end{theorem}
Theorem \ref{t11} relates coprime uncertainty balls to balls defined in the time-varying gap metric.
In particular, maximizing the uncertainty radius for coprime uncertainty results in maximizing
the uncertainty radius quantified in the TV gap metric. In the next section, we study the problem
of robust stabilization for time-varying normalized coprime factor perturbations and obtain operator
theoretic results, which generalize similar results in \cite{glover,smith,GS2} known to hold for LTI systems.
This is equivalent to maximizing the uncertainty radius in the TV gap metric by Theorem \ref{t11}, and the results obtained
here apply to uncertainty quantified in terms of the TV the gap metric as well.
\section{Robust Stabilization Under Coprime Factor Uncertainty}
\label{s2}
If $r_{opt}$ is the supremum over all $r$ such that $C$ stabilizes $B_s(P, r)$ or equivalently $B(P,r)$, then for
TV systems the following inequality holds (\cite{avraham}, p. 257)
\beqa
\inf_{n \geq 0} \left\{ \inf_{Q \in B_c(\ell^2,\; \ell^2)} \left\| \vect{U+MQ}{V + N Q} (I-P_n)
\right\| \right\}
\leq \frac{1}{r_{opt}}  \leq \inf_{Q \in  \Eu {B}_c (\ell^2, \ell^2)} \left\| \vect{U+MQ}{V + N Q} \right\|
\nonumber \\
\label{17bis}
\eeqa
For time-invariant systems these numbers are equal. In the same vein as \cite{smith}, we define $r_o$ as the right-hand
side of (\ref{17bis}).
%
\beqa
r_o^{-1} = \inf_{Q \in  \Eu {B}_c (\ell^2, \ell^2)} \left\| \vect{U}{V} + \vect{M}{N} Q \right\|
\label{18}
\eeqa
We will see that the infimum is achieved for some $Q_o \in  \Eu {B}_c (\ell^2, \ell^2)$ \cite{cdc04_tv,djouadi}.
\\
\\
In the sequel we are concerned with solving the optimization (\ref{18}) which reciprocal is
an upper bound for the optimal robustness radius. The solution proposed applies as well
to the lower bound in (\ref{17bis}) by a restriction to a subspace, that is, by restricting
the operators introduced in the sequel to the range of $(I-P_n)$, $(I-P_n)\ell^2$. The inequality (\ref{18}) allows
the estimation of $r_{opt}$ by computing a lower and upper bound.
\\
\\
The commutant lifting theorem (CLT) has been proposed by Sz.Nagy and Foias \cite{nagy68,nagy70}.  A time-varying version which
corresponds to nest or triangular algebras will be used and is discussed next. Following \cite{dav,ball} a nest ${\cal N}$ of
a Hilbert space ${\check{H}}$ is a family of closed subspaces of $\check{H}$ ordered by inclusion. The triangular or nest algebra
${\cal{T}} ({\cal{N}})$ is the set of all operators $T$ such that $TN \subseteq N$ for every element $N$ in ${\cal{N}}$. A representation of
${\cal{T}} ({\cal{N}})$ is an algebra homomorphism $h$ from ${\cal{T}} ({\cal{N}})$ into the algebra $\Eu B({\cal H})$ of bounded linear operators
on a Hilbert space ${\cal H}$. A representation is contractive if $\| h(A)\| \leq \|A\|$, for all $A\in {\cal{T}} ({\cal{N}})$. It is weak$^\star$ continuous if $h(A_i)$ converges to zero in the weak$^\star$ topology of $\Eu B({\cal H})$ whenever the net $\{A_i\}$ converges to zero in the weak$^\star$ topology of $\Eu B(\check{H})$. The representation $h$ is said to be unital if $h(I_{\check{H}}) =I_{{\cal H}}$, where $I_{\check{H}}$ is the identity operator on $\check{H}$, and $I_{{\cal H}}$ the identity operator on ${\cal H}$. The Sz. Nagy Theorem asserts that any such a representation $h$ has a $\Eu B(\check{H})$-dilation, that is, there exists a Hilbert space ${\cal K}$ containing ${\cal H}$, and a positive representation $H$ of $\Eu B(\check{H})$ such that $P_{{\cal H}} H(A) \mid_{{\cal H}}= h(A)$, where $P_{{\cal H}}$ is the orthogonal projection from ${\cal K}$ into ${\cal H}$ \cite{ball,dav}.
We now state the CLT for nest algebras.
\begin{theorem} \label{th1} \cite{ball,dav}
Let
\beqn
&h:&  {\cal{T}} ({\cal{N}}) \longmapsto {\Eu{B}}({\cal H}) \\
&h^\star:& {\cal{T}} ({\cal{N}}) \longmapsto {\Eu{B}}({\cal H}^{\prime} )
\eeqn
be two unital weak$^\star$ continuous contractive representations with $\Eu B (\check{H})$-dilations
\beqn
&H: & \Eu B(\check{H}) \longmapsto \Eu B({\cal K}) \\
&H^\star: & \Eu B(\check{H}) \longmapsto \Eu B ({\cal K}^{\prime})
\eeqn
respectively. Assume that $X:\; {\cal H} \longmapsto {\cal H}^{\prime}$ is a linear operator with $\|X \| \leq 1$, such that $Xh(A) = h^\prime (A)X$ for all $A \in {\cal{T}} ({\cal{N}})$, that is, $X$ intertwines $h$ and $h^{\prime}$. Then there exists an operator $Y:\; {\cal K} \longmapsto {\cal K}^{\prime}$ such that
\begin{description}
\item[i)] $\| Y \| \leq 1$.
\item[ii)] $Y$ intertwines $H$ and $H^\prime$, that is, $YH(A) = H^\prime (A)Y$ for all $A \in \Eu B(\check{H})$.
\item[iii)] $Y$ dilates $X$, that is, $Y:\; {\cal M} \longmapsto {\cal M}^\prime$, and
$P_{{\cal H}^\prime} Y \mid_{{\cal M}} = X P_{{\cal H}} \mid_{{\cal M}}$, where ${\cal H} = {\cal M} \ominus {\cal N}$ is
the orthogonal representation of ${\cal H}$ as the orthogonal difference of invariant subspaces for $H\mid_{{\cal{T}} ({\cal{N}})}$,
and similarly for ${\cal H}^\prime$.
\end{description}
\end{theorem}
Observe that (as in the time-invariant case) the operator
\beqa
\check{Z} :=  \left( \br{cc} M^\star & N^\star  \\ -\hN
 & \hM \er \right) \in \Eu B(\ell^2 \times \ell^2,\; \ell^2 \times \ell^2)  \label{19}
\eeqa
and the operator induced norm
\beqa
\left\| \vect{A}{B} \right\| &:=& \sup_{f \in \ell^2,\; \| f\|_2\leq 1} \left\|\vect{A}{B} f \right\|_2 \nonumber \\
&=& \sup_{f \in \ell^2,\; \| f\|_2\leq 1} \Bigl( \|Af\|_2^2 + \|Bf\|_2^2 \Bigl)^{\frac{1}{2}}, \;\;
\vect{A}{B} \in \EuBt  \label{20}
\eeqa
is unitarily invariant, that is, for any unitary operator $\tilde{U} \in \EuBt$, we have
\beqn
\left\| \tilde{U} \vect{A}{B} \right\| = \left\| \vect{A}{B} \right\|
\eeqn
In particular, for $\tilde{U}=\check{Z}$, we have
\beqa
\left\| \vect{U}{V} + \vect{M}{N} Q \right\| &=& \left\|\check{Z} \vect{U}{V} +
\vect{M}{N} Q \right\| \nonumber
\\
&=& \left\|\vect{R+Q}{I} \right\| = \bigl(1 + \|R+Q\|^2 \bigl)^{\frac{1}{2}} \nonumber \\
\label{21}
\eeqa
where $R:= M^\star U + N^\star V \in \Eu B(\ell^2, \; \ell^2)$. Next, we show that
\beqa
r_o = \bigl(1 + \|H_R\|^2 \bigl)^{-\frac{1}{2}}
\label{22}
\eeqa
where $H_R$ is the time-varying Hankel operator with symbol $R$.
Note that as in the LTI case $0 < r_o \leq 1$.
\\ \\
To define the operator $H_R$, we need some mathematical preliminaries. Let
\beqn
Q_n := I-P_n,\;\; {\rm for}\;\; n=-1,0,1, \cdots
\eeqn
where $P_n$ is the standard truncation operator, which is also an orthogonal projection that
sets all outputs after time $n$ to zero, and
\beqn
P_{-1} := 0\;\; {\rm and}\;\;\; P_\infty := I
\eeqn
Then $Q_n$ is a projection, and we associate to it the following nest
\beqn
\mathcal{N} := \{ Q_{n} {\ell}^2, \; n=-1, 0, 1, \cdots \}
\eeqn
The space of causal bounded linear operators $\Bc$ can be viewed as a triangular or
nest algebra, which leaves invariant every subspace $N \in \mathcal{N}$, that is,
$TN \subseteq N$, $\forall T\in \Bc$. In fact, $\Bc$ can be written as
\beqa
\Bc 
=\{ A \in \Eu{B}(\ell^2, \ell^2):\;  (I-Q_n) A Q_n = 0, \; \forall \; n\}
\label{23}
\eeqa
Now, call ${\cal{C}}_2$ the class of compact operators on $\ell^2$ called the Hilbert-Schmidt
or Schatten 2-class \cite{schatten,dav} under the norm,
\beqn
\| A \|_2 \; :=\; \Bigl(tr (A^\star A)\Bigl)^{\frac{1}{2}}
\eeqn
Define the space
\beqa
{\cal A}_2 := {\cal C}_2 \cap {\Eu{B}}_{c}({\ell^{2}}, {\ell^{2}})
\label{24}
\eeqa
Then ${\cal{A}}_2$ is the space of causal Hilbert-Schmidt operators. Operators in ${\Eu{B}}_{c}({\ell^{2}}, {\ell^{2}})$ may be considered as operators in ${\Eu{B}}_{c}({\cal{A}}_2, {\cal{A}}_2)$ and vise-versa \cite{ball}. Define the orthogonal projection ${\cal{P}}$ of ${\cal{C}}_2$ onto ${\cal{A}}_2$. ${\cal{P}}$ is the lower triangular truncation or nest projection. \\
\\
Let $\Pi_1$ be the orthogonal projection on the subspace
$({\cal A}_2 \oplus {\cal A}_2) \ominus \vect{M}{N} {\cal A}_2$ the
orthogonal complement of $\vect{M}{N} {\cal A}_2$ in the operator Hilbert
space ${\cal A}_2 \oplus {\cal A}_2$ under the inner product
\beqn
(A, \; B ) := tr (B^\star A),\;\; A,\; B \in {\cal A}_2 \oplus {\cal A}_2
\eeqn
In the following Lemma the orthogonal projection $\Pi_1$ is computed explicitly.
\begin{lemma}\label{lemma1}
\beqa
\Pi_1 = I-\vect{M}{N}{\cal{P}} \bigl(M^\star,\; N^\star\Bigl)
\label{orth}
\eeqa
\end{lemma}
{\bf Proof.}
Call $K:=\vect{M}{N}$. For $Z \in {\cal A}_2 \oplus {\cal A}_2$, let us compute
\beqa (I-K {\cal P}  K^\star )^2 Z
&=& ((I-K {\cal P}  K^\star )(I-K {\cal P}  K^\star )Z \nonumber \\ &=& (I-K{\cal P} K^\star -
K{\cal P}  K^\star +K{\cal P} K^\star {\cal P} K^\star )Z \nonumber \\
&=& (I-2K{\cal P} K^\star +K {\cal P} K^\star )Z \nonumber \\ && {\rm
since}\;\; K^\star K=I \;\;\; {\rm and }\;\;\; {\cal P}^2= {\cal P}
\nonumber \\ &=&(I-K{\cal P} K^\star ) Z \label{3.33} \eeqa so $(I-K{\cal P}
K^\star )$ is indeed a projection. \\ Clearly the adjoint $(I-K{\cal P}
K^\star )^\star$ of $(I-K{\cal P} K^\star )$ is equal to $(I-K{\cal P} K^\star )$ itself, so that
$(I-K{\cal P} K^\star )$ is an orthogonal projection. Next we show that the null space of $(I-K{\cal P} K^\star )$,
$Ker(I-K{\cal P} K^\star ) = K {\cal A}_2$. \\
Let $Z \in Ker(I-K{\cal P} K^\star )$ then 
$(I-K{\cal P} K^\star )Z =0 \Longrightarrow Z = K{\cal P} K^\star Z$, 
since $K^\star Z \in {\cal{C}}_2$, then ${\cal P} K^\star Z \in
{\cal A}_2$ and therefore $Z \in K {\cal A}_2$. Hence
$Ker(I-K{\cal P} K^\star ) \subset K {\cal A}_2$. Conversely, let $Z \in {\cal A}_2$,
then
\beqn (I-K{\cal P} K^\star )KZ = KZ - K{\cal P} Z =
KZ - KZ = 0
\eeqn
Thus $KZ \in Ker(I-K{\cal P} K^\star )$, so
$K {\cal A}_2 \subset Ker(I-K{\cal P} K^\star )$, and therefore
$(I-K{\cal P} K^\star ) = K {\cal A}_2$, and the Lemma is proved.
Following \cite{power} an operator $X$ in ${\cal B}(\ell^2, \ell^2)$ determines a Hankel
operator $H_X$ on ${\cal A}_2$ if
\beqa
H_X A = (I-{\cal P}) X A, \;\;\; {\rm for}\; A \in {\cal A}_2
\label{25}
\eeqa
For $R\in \EuBt$ the TV Hankel operator solves the minimization
\cite{cdc04_tv,djouadi}
\beqa
\inf_{Q\in \Bc} \| R + Q \| = \| H_R\|
\label{26}
\eeqa
where $H_R = (I-{\cal P}) (M^\star U + N^\star V)$. Formula (\ref{26}) is the time-varying analogue of the standard Nehari problem.
\\
The following expression for $r_o$ was obtained in the LTI case was obtained
in \cite{glover} using state space techniques, and in \cite{smith} using operator theory.
We will give its time-varying counterpart and give an operator theoretic proof along the lines of \cite{smith}.
\begin{theorem}
\label{t4}
\beqa
r_o^{-2} =  1- \Bigl\| H_{\vect{\hM^\star}{\hN^\star}} \Bigl\|^2
\label{27}
\eeqa
The quantity $r_o$ is a bound on the maximal achievable TV stability margin under coprime factor uncertainty
for LTV systems, and reduces to the maximal stability margin for LTI system.
\end{theorem}
{\bf Proof.}
Call $S := ({\cal A}_2 \oplus {\cal A}_2) \ominus \vect{M}{N} {\cal A}_2$,
and define the operator
\beqa
\Xi &:& {\cal A}_2 \longmapsto  S \nonumber \\ 
\Xi &:=& \Pi_1 \vect{U}{V} \label{28}
\eeqa
For the lower bound in (\ref{17bis}) it suffices to restrict the operator $\Xi$ to the subspace $(I-P_n){\cal A}_2$.
Using the CLT \cite{dav,ball,partington}, it follows
\beqa
r_o^{-1} = \| \Xi \|
\label{29}
\eeqa
To see this we need a representation of $\Bc$, that is, an algebra homomorphism, say, $h(\cdot)$
(respectively $h^\prime(\cdot)$), from $\Bc$, into the algebra
$B({\cal A}_2, {\cal A}_2)$ (respectively $\Eu B_c(S, S)$), of bounded linear operators from
${\cal A}_2$ into ${\cal A}_2$
$\bigl($respectively from $S$ into $S\bigl)$. Define the representations $h$ and $h^\prime$ by
\beqa
&h:& \Bc \longmapsto B({\cal A}_2, {\cal A}_2), \;\;
h^\prime: \Eu B_c(\ell^2, \ell^2) \longmapsto \Eu B_c(S, S)
\\
&h(A)& := R_A, \;\; A \in \Bc, \;\; h^\prime(A) := \Pi_1 R_A, \;\; A \in
\Eu B_c(\ell^2, \ell^2)  \nonumber
\eeqa
where $R_A$ denotes the right multiplication associated to the operator $A$ defined on the
specified Hilbert space. By the Sz. Nagy dilation Theorem there exist dilations $H$
(respectively $H^\prime$) for $h$ (respectively $h^\prime$) given by
\beqa
&H(A)& = R_A \;\; {\rm on} \;\; {\cal A}_2 \;\; {\rm for}\;\; A \in \Eu B(\ell^2, \ell^2)
\\
&H^\prime(A)& = R_A \;\; {\rm on} \;\;{\cal A}_2 \oplus {\cal A}_2   \;\; {\rm for}\;\; A \in
\Eu B (\ell^2, \ell^2)
\eeqa
The spaces ${\cal A}_2$ and $S$ can be written as orthogonal differences of subspaces invariant under
$H$ and $H^\prime$, respectively, as
\beqa
{\cal A}_2 = {\cal A}_2 \ominus \{0\}, \;\;\; S = {\cal A}_2 \oplus {\cal A}_2 \ominus \vect{M}{N} {\cal A}_2
\eeqa
Now we have to show that the operator $\Xi$ intertwines $h$ and $h^\prime$, that is, $h^\prime(A) \Xi = \Xi h(A)$ for all $A \in \Bc$,
\beqn
h^\prime(A) \Xi &=& \Pi_1 R_A \Pi_1 \vect{U}{V}\mid_{{\cal A}_2} = \Pi_1 R_A \vect{U}{V}\mid_{{\cal A}_2} \\
&=& \Pi_1 \vect{U}{V} R_A \mid_{{\cal A}_2} = \Xi h(A)
\eeqn
Applying Theorem \ref{th1} implies that $\Xi$ has a dilation $\Xi^\prime$ that intertwines $H$ and $H^\prime$, i.e., $\Xi^\prime H (A) = H^\prime (A) \Xi^\prime$, $\forall A \in B(\ell^2, \ell^2)$. By Lemma 4.4. in \cite{ball} $\Xi^\prime$ is a left multiplication operator acting from ${\cal A}_2$ into ${\cal A}_2 \oplus {\cal A}_2$, and causal. That is, $\Xi^\prime = L_K$ for some $K \in \Eu B_c({\cal A}_2, \;{\cal A}_2 \oplus {\cal A}_2)$, with $\| K\| = \|\Xi^\prime \|= \|\Xi \|$. Call $\tilde{U} := \vect{U}{V}$, then $\Xi = \Pi_1 \tilde{U} = \Pi_1 K$, which implies $\Pi_1 (\tilde{U}-K) = 0$. Hence, $ (\tilde{U}-K) f \in \tilde{M} {\cal A}_2$, where $\tilde{M}:=
\vect{M}{N}$ and for all $f \in {\cal A}_2$. That is, $ (\tilde{U}-K) f = \tilde{M}g,\; \exists g\in {\cal A}_2$, which can be written as $\tilde{M}^\star (\tilde{U}-K) f =g\in {\cal A}_2$.
In particular, $\tilde{M}^\star (\tilde{U}-K) f \in \Eu B_c(\ell^2, \ell^2)$, for all $f \in \Eu B_c(\ell^2, \ell^2)$ of finite rank. By Theorem 3.10 \cite {dav} there is a sequence $F_n$ of finite rank contractions in $\Eu B_c(\ell^2, \ell^2)$ which converges to the identity operator in the strong *-topology.  By an approximation argument it follows that
$\tilde{M}^\star (\tilde{U}-K) \in \Eu B_c(\ell^2, \; \ell^2)$. Letting $Q:= \tilde{M}^\star (\tilde{U}-K)$ we have
$g = Qf$. We conclude that $\tilde{U}-K =\tilde{M} Q$, that is, $\tilde{U}- \tilde{M} Q = K$, with $\|K\|=\|\Xi\|$ as required.
\\
\\
A consequence of the commutant lifting Theorem is that there exists an optimal $Q_o \in \Eu B_c(\ell^2, \; \ell^2)$ such that
the infimum in (\ref{18}) is achieved. Moreover, Lemma \ref{lemma1} implies that the operator $\Xi$ is given by the following
analytic expression
\beqa
\Xi =  \vect{U}{V} -\vect{M}{N}{\cal{P}}\bigl(M^\star U + N^\star V\bigl)
\eeqa
and the subspace $S$ is given by
\beqa
S = \Bigl({\cal A}_2 \oplus  {\cal A}_2 \Bigl)-\vect{M}{N}{\cal{P}} (M^\star,\; N^\star )  \Bigl({\cal A}_2 \oplus  {\cal A}_2 \Bigl)
\eeqa
Next, let $\Gamma$ be the left multiplication operator on $S$ as
\beqa
\Gamma &:& S \longmapsto {\cal A}_2 \nonumber \\
\Gamma &:=& (-\hN,\; \hM)|_S
\label{30}
\eeqa
Then, for $X \in {\cal A}_2$ we have
\beqn
\Xi X = \vect{U}{V} X - \vect{M}{N} Y \in S
\eeqn
for some $Y \in {\cal A}_2$ and
\beqn
\Gamma (\Xi X) = (-\hN,\; \hM) \left[ \vect{U}{V} X - \vect{M}{N} Y \right] = X
\eeqn
therefore $\Gamma\Xi = I$ on $S$. Moreover, for $s\in S$, we have
\beqn
\Gamma s &=& (-\hN,\; \hM)  s \\
&=& (-\hN,\; \hM) \left[ \vect{s_1}{s_2} - \vect{M}{N} s_3 \right]
\eeqn
for some $s_1, s_2$ and $s_3$ in ${\cal A}_2$. Now, apply $\Xi$ to get
\beqn
\Xi(\Gamma s) = \vect{U}{V} (-\hN,\; \hM)  s - \vect{M}{N} (-\hN,\; \hM)  s
\eeqn
which by (\ref{6}) is equal to
\beqn
\Xi(\Gamma s)=\vect{U}{V} (-\hN,\; \hM)  s - \vect{M}{N} (\hat{V},\; \hat{U})  s
\eeqn
and thus
\beqn
\Xi(\Gamma s)  = s
\eeqn
showing that $\Xi\Gamma =I$ on $S$, $\Gamma\Xi = I$ on ${\cal A}_2$, that is $\Gamma$
is the inverse of $\Xi$, $\Gamma = \Xi^{-1}$. This implies
\beqa
\|\Xi \| = \inf_{s\in S,\; \|s\|_2\leq 1} \|\Gamma s\|^{-1} =: \tau(\Gamma)^{-1}
\eeqa
Now, note that the adjoint operator of $\Gamma$, $\Gamma^\star$ is defined by
\beqa
\Gamma^\star &:& {\cal A}_2 \longmapsto S \label{31} \\
\Gamma^\star X &=&  \Pi_1 \vect{-\hN^\star}{\hM^\star}  X, \;\; X \in {\cal A}_2 \nonumber
\eeqa
and the TV Hankel operator $\Upsilon$ defined by
\beqn
\Upsilon &:& {\cal A}_2 \longmapsto \bigl( {\cal C}_2 \oplus {\cal C}_2\bigl) \ominus
\bigl( {\cal A}_2 \oplus  {\cal A}_2 \bigl) =: {\cal A}_2^{2\perp}
\nonumber \\
\Upsilon &=& (I-{\cal P}) \vect{-\hN^\star}{\hM^\star}
\eeqn
The subspace ${\cal{A}}_2^{2\perp}$ is the orthogonal complement of ${\cal{A}} \oplus {\cal{A}}$ in ${\cal C}_2 \oplus
{\cal C}_2$. Consider the left multiplication operator $\vect{-\hN^\star}{\hM^\star}$ defined from
${\cal A}_2$ into ${\cal C}_2$ which is isometric since we are using normalized coprime
factorizations. Then for any $X\in {\cal A}_2$ we have
\beqn
\vect{-\hN^\star}{\hM^\star} X &=& {\cal P} \vect{-\hN^\star}{\hM^\star} X+ (I-{\cal P})
\vect{-\hN^\star}{\hM^\star} X \\
&=&  \Pi_1 \vect{-\hN^\star}{\hM^\star} X + (I-{\cal P})
\vect{-\hN^\star}{\hM^\star} X \\
&=& \Gamma^\star X + \Upsilon X
\eeqn
But since $\Gamma \Gamma^\star + \Upsilon^\star \Upsilon =I$, that is, the operator
$\Gamma^\star + \Upsilon$ is an isometry, and $\tau(\Gamma^\star)^2 + \|\Upsilon\|^2 =1$.
Further, $\Xi^\star \Gamma^\star =I$ and so $\|\Gamma^\star\| \geq \frac{1}{\|\Xi\|}$,
i.e., $\Gamma^\star$ is bounded below implying $\|\tau(\Gamma)\|= \tau(\Gamma^\star)\|$. By definition
\beqa
\| \Upsilon \| = \Bigl\| H_{\vect{\hM^\star}{\hN^\star}} \Bigl\|
\eeqa
yielding
\beqa
 \Bigl\| H_{\vect{\hM^\star}{\hN^\star}} \Bigl\|^2 + \|\Xi\|^{-2} =1
\eeqa
which implies the result of the theorem (\ref{27}).
\\
{\bf Remark.} For LTI systems the connection between robust stabilization under coprime factor uncertainty and the Matrix-
Valued Corona problem has been pointed out in \cite{GS2}. It amounts to finding a left inverse of
$\vect{M}{N}$ with the norm bounded above by $\frac{1}{r}$. Robust stabilization for LTV systems involves an
Operator Corona Problem \cite{dav}. The latter may be stated in our case as follows: We seek causal
bounded linear operators $U,V \in \Eu {B}_c (\ell^2, \ell^2)$ such that
\beqa
\bigl(V,\; U\bigl)\vect{M}{N} = V{M}+ U{N}= I
\eeqa
A necessary and sufficient condition is given by the following Theorem.
\begin{theorem} \label{t3} The following statements are equivalent
\bi
\item[(i)] There exist operators $U,V \in \Eu {B}_c (\ell^2, \ell^2)$ such that
\beqa
V{M}+ U{N}= I\;\; {\rm and}\;\;\;  \left\| \vect{U}{V} \right\| \leq \frac{1}{r}
\eeqa
\item[(ii)]
\beqa
\sup_{n,\; f \in \ell^2} \frac{\|P_n f\|_2}{\|P_n M f +P_n N f\|_2} \leq \frac{1}{r} < \infty
\eeqa
\ei
\end{theorem}
{\bf Proof.} The Theorem follows from Theorem 4 in \cite{katsoulis}.
\\
\\
If the Hankel operator $\Upsilon$ 
is compact, then its induced norm 
is equal to its maximal singular value, and in fact, its norm is achieved by some operator $X_o \in {\cal{A}}_{2},\;
\|X_o\|_2 \leq 1$, i.e..
\beqa
\sup_{\|X\|_2 \leq 1} \| \Upsilon X\|_2 = \|\Upsilon X_o\|_2 = \|\Upsilon \|\|X_o\|_2 = \|\Upsilon \|
\eeqa
A necessary and sufficient condition for any Hankel operator $H_X$ to be compact is
that $X$ belongs to a triangular algebra $\Eu{A}$ plus the space of compact operator
denoted $\Eu{K}$, that is, $X \in {\Eu{A}} + \Eu{K}$ \cite{power}. In our case, if we define
$\Eu{K}$ as the space of compact operator from ${\cal{A}}_2$ into ${\Eu{C}}_2$ and
${\Eu{A}} := {\Eu{B}}({\cal{A}}_2, {\cal{A}}_2 \oplus {\cal{A}}_2)$, then
$\Upsilon$ is compact if and only if
\beqn
\vect{-\hN^\star}{\hM^\star} \in \Eu{A} + \Eu{K}
\eeqn
As in the LTI case \cite{GS2}, this has implications for the operators $\Xi$
and $\Gamma$, in particular $\Xi$ attains its norm. This is summarized in the following theorem.
\begin{theorem}  \label{t65}
Let $0 < \lambda < 1$. Then the following are equivalent
\bi
\item[1)] $\lambda$ is a singular value of $\Upsilon$.
\item[2)] There exists $X\neq 0 \in {\cal{A}}_2$, $Y^\star \neq 0  \in {\cal{A}}_2^{2\perp}$, $W \in
\Pi_1 ({\cal{A}}_2 \oplus {\cal{A}}_2)$ such that
\beqa
\vect{-\hN^\star}{\hM^\star} X - (1-\lambda^2)^\frac{1}{2} W = \lambda Y^\star, \;
(-\hN, \; \hM) Y^\star = \lambda X \label{32}
\eeqa
\item[3)] $(1-\lambda^2)^\frac{1}{2}$ is a singular value of $\Gamma$.
\item[4)] There exists $X\neq 0 \in {\cal{A}}_2$, $Y^\star \neq 0  \in {\cal{A}}_2^{2\perp}$, $W \in
\Pi_1 ({\cal{A}}_2 \oplus {\cal{A}}_2)$ such that
\beqa
(-\hN, \; \hM) W  &=& (1-\lambda^2)^\frac{1}{2} X \label{33} \\
\vect{-\hN^\star}{\hM^\star} X - \lambda Y^\star & = & (1-\lambda^2)^\frac{1}{2} W \label{34}
\eeqa
\item[5)] $(1-\lambda^2)^{-\frac{1}{2}}$ is a singular value of $\Xi$.
\ei
\end{theorem}
{\bf Proof.} Follows as in the LTI case \cite{GS2} and is omitted.\\\\
In other words, Theorem \ref{t65} says that $X$ and $Y^\star$ are the Schmidt pairs for the operator
$\Upsilon$, i.e.,
\beqa
\Upsilon X = \lambda Y^\star, \;\;
\Upsilon^\star Y^\star = \lambda X
\eeqa
Similarly,
\beqa
\Gamma W = (1-\lambda^2)^{\frac{1}{2}} X, \;\;
\Gamma^\star X = (1-\lambda^2)^{\frac{1}{2}} W
\eeqa
and therefore,
\beqa
\Xi X = (1-\lambda^2)^{-\frac{1}{2}} W, \;\;
\Xi^\star W = (1-\lambda^2)^{-\frac{1}{2}} X
\eeqa
The last expression will be used to give an explicit formula for the optimal solution.
In this section we show that the optimal $Q_o$ in (\ref{18}) exists and that the infimum is indeed achieved. To do so we
invoke duality theory \cite{djouadi}. We need first to introduce, the class of compact operators on $\ell^2$ called the
trace-class of operators acting from $\ell^2$ into $\ell^2\times\ell^2$, denoted ${\cal C}_1(\ell^2, \ell^2 \times \ell^2)$,
under the trace-class norm \cite{schatten,dav},
The following result gives a necessary and sufficient condition for $\lambda$ to be a singular
value of $\Upsilon$ and generalizes its LTI counterpart.
\begin{theorem}\label{t7}
Let $0 < \lambda < 1$. Then $\lambda$ is a singular value of $\Upsilon$ if
and only if there exists $W \neq 0 \in {\cal{A}}_2\oplus {\cal{A}}_2$ such that
\beqa
Z:= \vect{-\hN^\star}{\hM^\star} (-\hN, \; \hM) W - (1-\lambda^2)^\frac{1}{2}  W \in
{\cal{A}}_2^{2\perp} \label{34}
\eeqa
Moreover, if (\ref{34}) holds then $W \in \Pi_1 ({\cal{A}}_2 \oplus {\cal{A}}_2)$ and
Theorem \ref{t65}, 3) and 4) hold for
\beqa
X = \frac{1}{(1-\lambda^2)^\frac{1}{2}}(-\hN, \; \hM) W, \;\;
Y^\star = \frac{1}{\lambda(1-\lambda^2)^\frac{1}{2}}Z \label{35}
\label{36}
\eeqa
\end{theorem}
The optimal operator $Q_o \in \Bc$ which achieves the infimum in (\ref{18}) satisfies
\beqa
\left\| \vect{U}{V} + \vect{M}{N} Q_o \right\| =
\inf_{Q \in  \Eu {B}_c (\ell^2, \ell^2)} \left\| \vect{U}{V} + \vect{M}{N} Q \right\|
= \|\Xi\| \nonumber \eeqa
and can be computed from the identities
\beqa
\Xi  X_o = \|\Xi\|  W_o,\;\;
\Xi^\star W_o = \|\Xi\| X_o
\eeqa
that is, $X_o$ and $W_o$ is the Schmidt pair corresponding to the maximum singular value $\lambda=\|\Xi\|$,
and satisfies Theorems \ref{t65} and \ref{t7}. And the following operator identity gives $Q_o$,
\beqa
 \vect{U}{V}X_o + \vect{M}{N} Q_o X_o = \Xi X_o = \|\Xi \| W_o \\
 \Longleftrightarrow Q_o X_o = -\bigl( M^\star U + N^\star V \bigl) X_o + \|\Xi \|
 \bigl( M^\star, \; N^\star \bigl) W_o
\eeqa
which solution is guaranteed to exist.
\section{Conclusion} \label{s4}
In this paper we considered the problem of robust stabilization of LTV systems in the gap metric and its connection
with coprime factor uncertainty. We studied the problem of computing the optimal controller and characterize the
radius of the maximal TV gap ball about the plant which can
be stabilized by a fixed controller. An Operator Corona Theorem which gives a necessary and sufficient condition
for the existence of a robustly stabilizing controller was pointed out. We introduced TV Hankel and some related operators along with their singular values
and vectors that play a central role in optimizing the TV gap and in the computation of the robust stabilizing
controller. The maximal stability margin under coprime factor uncertainty for LTV systems is characterized as the norm of
these operators. Our results generalize similar results obtained for the case of LTI systems in \cite{GS2}, and in fact reduce
to them in that case.

\end{document}